\documentclass[letterpaper]{article}

\title{An extension of the Moser-Tardos algorithmic local lemma}
\author{Wesley Pegden\footnote{
Courant Institute of Mathematical Sciences, New York University,
251 Mercer St, Rm 921,
New York, NY 10012
Email: pegden@math.nyu.edu.  Partially supported by NSF MSPRF grant 1004696.}
}

\date{March 12, 2011}

\usepackage{amsmath,amsthm,amssymb,amscd,latexsym,pstricks,pst-node,pst-tree}
\usepackage{url}
\usepackage{cite}
\usepackage{mathtools}
\usepackage{algpseudocode}

\newcommand{\p}{\textrm{P}}
\newcommand{\e}{\textrm{E}}

\newcommand{\vbl}{\textrm{vbl}}

\newcommand{\pp}{{\rm P}}

\newcommand{\ppp}{{\cal P}}
\newcommand{\aaa}{{\cal A}}
\newcommand{\ccc}{{\cal C}}

\newcommand{\ttt}{{\cal T}}
\newcommand{\vvv}{{\cal V}}

\newcommand{\sbs}{\subset}

\newcommand{\bfrac}[2]{\frac{\displaystyle #1}{\displaystyle #2}}

\newtheorem{theorem}{Theorem}[section]

\newtheorem{lemma}[theorem]{Lemma}

\newcommand{\comments}[1]{}

{
\theoremstyle{definition}

}

{
\theoremstyle{remark}

}

\newcommand{\stm}{\setminus}

\begin{document}
\maketitle

\begin{abstract}
  A recent theorem of Bissacot, et al.~proved using results about the cluster expansion in statistical mechanics extends the Lov\'asz Local Lemma by weakening the conditions under which its conclusions holds.  In this note, we prove an algorithmic analog of this result, extending Moser and Tardos's recent algorithmic Local Lemma, and providing an alternative proof of the theorem of Bissacot, et al.~applicable in the Moser-Tardos algorithmic framework.
\end{abstract}

\section{Introduction}
If events $A_1,A_2,\dots,A_n$ are independent, then we have 
$\pp(\bigcap \bar A_i))>0$
so long as $\pp(A_i)<1$ for each $i$.  A central tool in probabilistic combinatorics is the Lov\'asz Local Lemma proved by Erd\H{o}s and Lov\'asz\cite{EL}, which can be seen as generalizing this simple fact to situations where some dependencies among the $A_i$ are allowed, in exchange for better bounds on the probabilities $\pp(A_i)$.

The Local Lemma is commonly presented through the framework of a \emph{dependency graph} on the events $A_i$, where if $\ccc$ is any family of non-neighbors of some $A_i$, then we have that $A_i$ is independent of the family $\ccc$ of events.  The Lov\'asz Local Lemma is then as follows:

\begin{theorem}[Lov\'asz Local Lemma]
  Let $G$ be any dependency graph for a finite family $\aaa$ of events, and suppose that there are real numbers $0<x_A<1$ $(A\in \aaa)$  such that for all $A\in \aaa$ we have
  \begin{equation}
    \pp(A)\leq x_{A}\prod_{B\sim A}(1-x_B).
\label{c.L}
  \end{equation}
Then 
\[
\pp\left(\bigcap_{A\in \aaa}\bar A\right)>\prod_{A\in \aaa}(1-x_A),
\]
and so in particular, we have
\begin{equation}
\pp\left(\bigcap_{A\in \aaa}\bar A\right)>0.
\end{equation}
\end{theorem}

The first breakthrough in finding an algorithmic version of the Local Lemma was made by Beck, who demonstrated his method on the classical Local Lemma application to 2-colorable hypergraphs.  Beck's method was subsequently refined and given a more general framework \cite{Alimp,MR,CS,S}, but required stronger bounds on the probabilities of the events than were required by the nonalgorithmic version.

In Moser and Tardos' recent breakthrough paper\cite{MT}, they give an algorithmic proof of the Lov\'asz Local Lemma in a setting which is general enough for nearly all applications of the Lemma in combinatorics, with bounds identical to those required by the nonalgorithmic version.  In the framework Moser and Tardos consider, the events in $\aaa$ depend on some underlying set $\vvv$ of independent random variables, and they denote by $\vbl(A)$ $(A\in \aaa)$ the minimal set of random variables from $\vvv$ on which each $A$ depends; $A$ is said to be `violated' with respect to a particular evalation of the variables in $\vbl(A)$ if the event occurs for that evaluation.  A Moser-Tardos dependency graph is one which implies that if events $A$ and $B$ are nonadjacent, then $\vbl(A)$ is disjoint from $\vbl(B)$.  (Note that this notion of a dependency graph is more restrictive than the Lov\'asz version based on probabilistic independence, as is demonstrated by an example of Kolipaka and Szegedy\cite{KS}.)   Moser and Tardos's theorem is then the following:
\begin{theorem}[Moser and Tardos]
  Let $\vvv$ be a finite set of mutually independent variables in a probability space, and let $\aaa$ be a finite family of events determined by these variables.  If there are real numbers $0<x_A<1$ ($A\in \aaa$) such that
\begin{equation}
\pp(A)\leq x_A\prod_{B\sim A}(1-x_B)
\label{c.MT}
\end{equation}
then there exists an assignment to the variables $\vvv$ which corresponds to no occurrence of any event from $\aaa$.  Moreover, the randomized algorithm described below resamples an event $A$ at most an expected $\frac{x_A}{1-x_A}$ times before finding the evaluation, thus the total number of resampling steps is $\sum_{A\in \aaa}\frac{x_A}{1-x_A}$ in expectation.
\end{theorem}
\noindent The Moser-Tardos algorithm consists just of beginning with a random evaluation of all the variables in $\vvv$, and then resampling $\vbl(A)$ for any violated events $A$ until no violated events remain.  Of course, the efficiency of the algorithm depends on the ability to resample variables efficiently and check whether individual events are violated; this is generally an easy implementation problem, however, making the analysis of the number of resampling steps the important issue.

Recently, Bissacot, Fern\'andez, Procacci and Scoppola proved the following improvement of the Lov\'asz Local Lemma:
\begin{theorem}[Bissacot, et al.\cite{BFPS}]
\label{t.ball}
  Consider a finite family $\aaa$ of events in some probability space $\Omega$, with some dependency graph $G$.   If there are real numbers $0<\mu_A<\infty$ such that
\begin{equation}
\p(A)\leq 
\bfrac{\mu_A}{\sum_{\substack{I\sbs \bar\Gamma(A)\\ I \textrm{ \emph{indep.}}}} \prod_{B\in I}\mu_B}
\label{c.ball}
\end{equation}
then $\p\left(\bigcap\limits_{A\in \aaa}\bar A\right)>0$.
\end{theorem}
\noindent It is not difficult to check that condition (\ref{c.ball}) is weaker than condition (\ref{c.L}) by considering the substitution $\mu_A=\frac{x_A}{1-x_A}$.  (Condition (\ref{c.L}) would be equivalent to (\ref{c.ball}) without the condition in the sum that the sets $I$ be independent.)  In \cite{BFPS}, they also give examples where this theorem improves some classical theorems proved with the Local Lemma.  Theorem \ref{t.ball} has also since been applied to improve some theorems on graph colorings in \cite{nlapp}.

Their proof of Theorem \ref{t.ball} is based on Shearer's characterization of labeled dependency graphs to which the conclusion of the Local Lemma applies\cite{Sh} and two of those authors' recent results on the radius of convergence of logs of partition functions\cite{cluster}.  (The connection between the Local Lemma and the partition functions of statistical mechanics was first made by Scott and Sokal\cite{SS}.)

In this short note, we prove an algorithmic analog to the result of Bissacot, et.~al.  That is, we will show that in the setting of Moser and Tardos's algorithmic Local Lemma, Moser and Tardos's bounds on the running time of their algorithm hold even with their condition (\ref{c.MT}) replaced with condition (\ref{c.ball}):

\begin{theorem}
\label{t.algimp}
  Let $\vvv$ be a finite set of mutually independent variables in a probability space, and let $\aaa$ be a finite family of events determined by these variables.  If there are real numbers $0<\mu_A<\infty$ ($A\in \aaa$) such that
\begin{equation}
\p(A)\leq 
\bfrac{\mu_A}{\sum_{\substack{I\sbs \bar\Gamma(A)\\ I \textrm{ \emph{indep.}}}} \prod_{B\in I}\mu_B},
\label{l.icond}
\end{equation}
then there exists an assignment to the variables $\vvv$ which corresponds to no occurrence of any event from $\aaa$.  Moreover, the Moser-Tardos algorithm resamples an event $A$ at most an expected $\mu_A$ times before finding the evaluation, thus the total number of resampling steps is $\sum\limits_{A\in \aaa}\mu_A$ in expectation.
\end{theorem}
\noindent (The running time bound here is equivalent to the Moser-Tardos bound under the substituation $\mu_A=\frac{x_A}{1-x_A}$.)

The proof of Theorem \ref{t.algimp} consists simply of re-doing one part of the proof of Moser and Tardos's theorem, taking advantage of some constraints which were not necessary for Moser and Tardos's original result.

Theorem \ref{t.algimp} can be seen as doing two things: first, it extends the result of Moser and Tardos by giving a weaker condition under which the identical result holds---note that this has also been done in a more general sense by Kolipaka and Szegedy\cite{KS}, who directly connect Shearer's condition with the Moser/Tardos algorithmic framework. Secondly, it gives an alternative proof of the result of Bissacot, et al.~(in the slightly more restrictive algorithmic setting) which is independent of Shearer's theorem and the cluster expansion methods used in \cite{cluster}.

\bigskip 

Bissacot, et al.~note that their Theorem \ref{t.ball} can be extended to \emph{Lopsided dependency graphs}, first considered by Erd\H{o}s and Spencer in \cite{ES}.  In their paper on their algorithmic Local Lemma, Moser and Tardos define an analog of lopsidependency in the algorithm/variable setting, and a reader familiar with Moser and Tardos's paper can easily verify that our improvement to Moser and Tardos's theorem applies to their theorem on algorithmic lopsided dependency graphs as well, as we only re-do their branching argument, which is applied to the lopsided case in the same way as in their main result.

\section{Proof}
\label{s.proof}

The Moser-Tardos algorithm is as follows:
\begin{algorithmic}[1]
 \Procedure{Moser-Tardos}{$\ppp$}
  \For {all $P\in \ppp$}
  \State $v_P\gets $(random evaluation of $P$)
  \EndFor
  \While {$\exists A$ s.t. $A$ is violated when $P=v_P$ ($\forall P$)} \label{l.w}
   \For {all $P\in \vbl(A)$}
   \State $v_P\gets $(new random evaluation of $P$)
  \EndFor
   \EndWhile
  \EndProcedure
\end{algorithmic}
\noindent (Note that when multiple events exist satisfying line \ref{l.w}, one of the satisfying events is chosen arbitrarily.)

Moser and Tardos' proof that this algorithm terminates in polynomial time (under condition (\ref{c.MT})) is based on the notion of a `witness tree'.  As the algorithm runs and bad events are found and resampled, a witness tree is assigned to each step of the algorithm (where a step consists of a resampling of an event).  A witness tree is a rooted tree with labels from $\aaa$.  The witness tree $W_t$ for step $t$ of the algorithm is constructed as follows: choose as its root a vertex labeled with whatever event $A_0$ was resampled at step $t$.  If the event $A_1$ which was resampled at step $t-1$ overlaps the label of the root, a vertex is added as a child of the root labeled with $A_1$. (We may have $A_1=A_0$.)  In general, for each step $i=t-1,t-2,\dots,1$ of the algorithm, if the event $A_i$ which was added at step $i$ overlaps any of the events currently labels of vertices of our partially constructed $W_t$, we add a vertex labeled with $A_i$ as the child of a vertex of maximum depth whose label overlaps $A_i$.  In the result, $W_t$, children always overlap their parents, and children of a common parent always get distinct labels (otherwise, whichever was added after the other would have been added as a child of the other).  Any tree $T$ with labels from $\aaa$ with these two properties is called a \emph{proper witness tree}.

Moser and Tardos's proof of their algorithm's efficiency consists of two parts: first, they show that any proper witness tree $T$ has probability at most
\begin{equation}
\prod_{v\in T} \p(A_v)
\end{equation}
of occurring as a witness tree at any point in the running of the algorithm, where here $A_v$ denotes the event labeling the vertex $v$.

Now, if an event $A$ is resampled at step $t$, the number of occurrences of $A$ as a label in the witness tree $W_t$ is equal to the number of times $A$ has been resampled on steps $1,\dots,t$---in particular, all witness trees which will occur in a run of the algorithm will be distinct.  Thus if we let $\ttt_A$ denote the set of proper witness trees with root label $A$, the expected value of the number $N_A$ of resamplings of $A$ which occur in a run of the algorithm is equal to
\begin{equation}
\e(N_A)=\sum_{T\in \ttt_A} \p(T \textrm{ occurs in the log})\leq \sum_{T\in \ttt_A}\prod_{v\in T} \p(A_v).
\label{l.nasum}
\end{equation}

The second part of Moser and Tardos's proof consists of bounding the sum of products in line (\ref{l.nasum}).  They do this by considering a random process for constructing trees:  Suppose $x_A$ ($A\in \aaa$) are real numbers between 0 and 1.  Fix now any event $A_0$.  In the first round of the process, a vertex labeled $A_0$ is created.  In each subsequent round, for each event vertex $v$ with label $A_v$ created in the previous round, and for each event $A_u\in \bar\Gamma(A_v)$ (in the dependency graph), a vertex $u$ with label $A_u$ is added as a child of $v$ with probability $x_{A_u}$.  (All of these choices are are made independently.)

Moser and Tardos prove:
\begin{lemma}[Moser Tardos Branching Lemma]
  For any proper witness tree $T$ with root labeled $A_0$, the probability $p_T$ that the process above produces exactly the tree $T$ is
\begin{equation}
  p_T=\frac{1-x_{A_0}}{x_{A_0}}\prod_{v\in T}\left(x_{A_v}\prod_{B\sim A_v}(1-x_B)\right).
\end{equation}
\qed
\end{lemma}

Thus, the Lemma gives us that 
\begin{equation}
1\geq \sum_{T\in \ttt_A} p_T\geq
\frac{1-x_{A}}{x_{A}}\sum_{T\in \ttt_A}\prod_{v\in T}\left(x_{A_v}\prod_{B\sim A_v}(1-x_B)\right)
\end{equation}
Thus the bound $P(A)\leq \left(x_{A}\prod_{B\sim A}(1-x_B)\right)$ for all $A$ implies, together with line (\ref{l.nasum}), that 
\begin{equation}
  \e(N_A)\leq \frac{x_A}{1-x_A}.
\end{equation}

Our improvement comes just from a slightly more careful branching argument.  Note that any witness tree which occurs in the log of the algorithm has the property that any children of a common vertex have labels which are nonadjacent in the dependency graph.  This condition---let's call it \emph{strongly proper}---is stronger than requiring simply that children be distinct.  Thus, we can strengthen line \ref{l.nasum}, as we have the bound

\begin{equation}
\e(N_A)=\sum_{T\in \ttt^S_A} \p(T \textrm{ occurs in the log})\leq \sum_{T\in \ttt^S_A}\prod_{v\in T} \p(A_v).
\label{l.rnasum}
\end{equation}
where $\ttt_A^S\sbs \ttt_A$ is the set of strongly proper witness trees.  

To bound the sum in (\ref{l.rnasum}), we consider a modified branching process which proceeds as follows.

Given real numbers $0<\mu_A<\infty$, we define $x_A=\frac{\mu_A}{\mu_A+1}$ (note that $0<x_A<1$) and fix any event $A_0$.  In the first round of the process, a vertex labeled $A_0$ is created.  In each subsequent round, for each event vertex $v$ with label $A_v$ in the previous round, we carry out a `subprocess', where for each $A_u\in \bar\Gamma(v)$ (in the dependency graph), a vertex $u$ with label $A_u$ is added as a child of $v$ with probability $x_{A_u}$ (the choices are independent).  At the end of the subprocess, we check if the label-set of the resulting set of children for $v$ is an independent set in the dependency graph.  If it is not, we delete the children created and restart the subprocess.  Note that $x_A<1$ (for all $A$) implies that the subprocess will eventually end (with probability 1) having produced an independent set.

Note that the process described above is equivalent to one in which, in each round and for each vertex $v$ from the previous round, we create a set of children $u$ with labels from a set chosen from all independent sets $I_v\sbs \bar\Gamma(v)$, where the likelihood of the choice of each independent set $I_v$ is weighted according the the product 
\[
w(I_v)=\left(\prod_{u\in I_v}x_{A_u}\right)\left(\prod_{u\in \bar\Gamma(v)\stm I_v}(1-x_{A_u})\right).
\]

\begin{lemma}[Improved Branching Lemma]
  For any strongly proper witness tree $T$ with root labeled $A_0$, the probability $p'_T$ that the modified branching process described above produces exactly the tree $T$ is
\begin{equation}
p'_T=\mu_{A_0}^{-1}
\prod_{v\in T}\bfrac{\mu_{A_u}}{\sum_{\substack{I\sbs \bar\Gamma(A_v)\\ I \textrm{ \emph{indep.}}}} \prod_{A\in I}\mu_A}.
\label{l.ipbound}
\end{equation}
\end{lemma}
\begin{proof}
Letting $W_v=\bar\Gamma_G(A_v)\stm \ell(\Gamma^+_T(v))$, where $\ell(v)$ is the label of vertex $v$, we have
\[
p'_T=\prod_{v\in T}
\bfrac{
 \prod_{u\in \Gamma^+_T(v)}x_{A_u} \prod_{B\in W_v}(1-x_B)
}
{
\sum_{\substack{I\sbs \bar\Gamma(A_v)\\ I \textrm{ indep.}}} \prod_{A\in I}x_A\prod_{B\in \bar\Gamma_G(A_v)\stm I}(1-x_B)
}.
\]
This can be rewritten as
\[
p'_T=\prod_{v\in T}
\bfrac{
\prod_{u\in \Gamma^+_T(v)}\frac{x_{A_u}}{1-x_{A_u}}
}
{
\sum_{\substack{I\sbs \bar\Gamma(A_v)\\ I \textrm{ indep.}}} \prod_{A\in I}\frac{x_A}{1-x_A}
}
\]
by dividing the top and bottom by $\prod_{B\in \bar\Gamma_G(A_v)}(1-x_B)$.  Since taking the double product $\prod_{v\in T} \prod_{u\in \Gamma^+_T(v)}$ is a equivalent to taking a product $\prod_{v\in T\stm \{v_0\}}$, where $v_0$ denotes the root vertex of $T$, this gives line (\ref{l.ipbound}), recalling that $x_A=\frac{\mu_A}{\mu_A+1}$ and so $\frac{x_A}{1-x_A}=\mu_A$.
\end{proof}
This is applied now in the same way as the branching lemma used by Moser and Tardos, but with regards to the family $\ttt^S_A$ of strongly proper witness trees rooted with $A$, instead of the family $\ttt_A$ of proper witness trees rooted with $A$.
We have
\begin{equation}
1\geq \sum_{T\in \ttt^S_A}p'_T\geq 
\mu_{A_0}^{-1}\sum_{T\in \ttt^S_A}
\prod_{v\in T}\bfrac{\mu_{A_u}}{\sum_{\substack{I\sbs \bar\Gamma(A_v)\\ I \textrm{ indep.}}} \prod_{A\in I}\mu_A}.
\label{}
\end{equation}
Putting this together with line (\ref{l.rnasum}), we see then that the condition (\ref{l.icond}) of the theorem implies that
\begin{equation}
  \e(N_A)\leq \mu_A,
\end{equation}
completing the proof the the Moser-Tardos algorithm still terminates in expected time
\[
\sum_{A\in \aaa} \mu_A.
\]\qed

\subsubsection*{Acknowledgement}
I'd like to thank Joel Spencer for some helpful discussions on this note.

\bibliographystyle{abbrv}

\end{document}